\documentclass{article}

\usepackage[latin1]{inputenc}
\usepackage[english]{babel}
\usepackage{amsmath,amssymb}
\usepackage[dvips]{graphicx}
\usepackage{fancyhdr}
\usepackage{enumerate}

\newtheorem{Lem}{Lemma}
\newenvironment{lem}[1][]{\begin{Lem}\begin{normalfont}\emph{#1 }}
{\end{normalfont}\finenu\end{Lem}}
\newtheorem{Teo}[Lem]{Theorem}
\newenvironment{teo}[1][]{\begin{Teo}\begin{normalfont}\emph{#1 }}
{\end{normalfont}\finenu\end{Teo}}
\newtheorem{Cor}[Lem]{Corollary}
\newenvironment{cor}[1][]{\begin{Cor}\begin{normalfont}\emph{#1 }}
{\end{normalfont}\finenu\end{Cor}}
\newtheorem{Pro}[Lem]{Proposition}
\newenvironment{pro}[1][]{\begin{Pro}\begin{normalfont}\emph{#1 }}
{\end{normalfont}\finenu\end{Pro}}
\newtheorem{Defi}[Lem]{Definition}
\newenvironment{defi}[1][]{\begin{Defi}\begin{normalfont}\emph{#1 }}
{\end{normalfont}\fine\end{Defi}}

\newenvironment{myproof}[1][]{\par\noindent\textbf{Proof{#1}. }}{\finedim\par}

\newcommand{\fine}{}
\newcommand{\finenu}{}
\newcommand{\finedim}{\hfill$\blacksquare$}

\newcommand{\Banach}{{\mathbb{R}^d}}
\newcommand{\dBanach}{{\mathfrak{X}^{*}}}
\newcommand{\Closed}[2][]{\mathbb{K}_{#1}}
\newcommand{\Fuzzy}[2][]{\mathbb{F}_{#1}}
\newcommand{\salgebra}{\mathfrak{F}}
\newcommand{\borel}[1][]{\mathcal{B}_{#1}}
\newcommand{\prob}[1][]{\mathbb{P}_{#1}}
\newcommand{\misura}{\mu}
\newcommand{\dd}{{\rm d}}

\newcommand{\unitSphere}{S^{\, d-1}}
\newcommand{\support}[2][]{h_{#2} #1 }

\newcommand{\Mink}{\oplus}
\newcommand{\Minkowski}{+}
\newcommand{\HausDist}{\delta_H}
\newcommand{\indicator}[1]{\mathbb{I}_{#1}}
\newcommand{\frv}{{FRV}}

\begin{document}

\title{A Note on Fuzzy Set--Valued Brownian Motion}
\author{Enea G. Bongiorno (enea.bongiorno@unimi.it)}
\maketitle

\abstract{In this paper, we prove that a fuzzy set--valued Brownian motion $B_t$, as defined in \cite{li:gua07}, can be handle by an $\Banach$--valued Wiener process $b_t$, in the sense that $B_t =\indicator{b_t}$; i.e. it is actually the indicator function of a Wiener process.}

\section{Introduction}

Stochastic (fuzzy) set--valued evolution is a relevant topic that was studied largely by different authors (e.g. \cite{li:gua07,li:ogu:kre02,mol05} and references therein). The following question was stated by Molchanov in \cite[Open Problem 1.24, p.316]{mol05}:
\begin{quote}
Define a set--valued analogue of the Wiener process and the corresponding stochastic integral.
\end{quote}
In \cite{li:gua07}, the authors tackle the proposed problem defining a fuzzy set--valued Brownian motion in $\Fuzzy[kc]{\Banach}$, the family of convex fuzzy subsets of $\Banach$ with compact support. In the sequel we shall prove that such a process is equivalent to consider simply a Wiener process in $\Banach$. This is based upon the fact that the Brownian motion is a zero--mean Gaussian (fuzzy set--valued) process.
\\
In fact, it is widely known (cf. \cite[Theorem~6.1.7]{li:ogu:kre02}) that a Gaussian random fuzzy set decomposes according to
\begin{equation}\label{eq:gauss_decomposition}
X=\mathbb{E} X \Mink \indicator{\xi },
\end{equation}
where $\mathbb{E} X$ is in the Aumann sense, $\xi$ is a Gaussian random element in $\Banach$ with $\mathbb{E} \xi=0$ and $\indicator{A}:\Banach\to\{0,1\}$ denotes the indicator function of any $A\subseteq \Banach$
\[
\indicator{A}(x)=\left\{
\begin{array}{ll}
1, & \textrm{if } x\in A,\\
0, & \textrm{otherwise},
\end{array}
\right.
\]
(for the sake of simplicity, whenever $A=\{a\}$ is a singleton we shall write $\indicator{ a }$ instead of $\indicator{\{a\} }$). Equation~\eqref{eq:gauss_decomposition} means that $X$ is just its expected value $\mathbb{E} X$ up to a random Gaussian translation $\xi$. In some sense, $\mathbb{E} X$ represents the ``deterministic'' part of $X$ whilst $\xi$ represents its random part.
It is also known (cf. \cite[Proposition 1.30, p.161]{mol05}) that a zero--mean random set is actually a random element in $\Banach$ with zero--mean. Such a result can be easily extended to the fuzzy case and, jointly to decomposition \eqref{eq:gauss_decomposition}, implies
\[
X=\indicator{0 } \Mink \indicator{\xi } = \indicator{\xi }.
\]
Roughly speaking, the definition of Brownian motion in \cite{li:gua07} for random fuzzy sets drives down the complexity of the chosen (fuzzy) framework. In fact, a Gaussian fuzzy random set with zero--mean is reduced to be a random Gaussian element in $\Banach$.
\\
In this paper we shall provide an alternative proof of the last fact using selections.

The paper is organized as follow. Section~\ref{sec: preliminary results} is devoted to preliminaries such as random (fuzzy) sets, embedding theorems and Brownian motion for fuzzy sets (according to \cite{li:gua07}). In Section~\ref{sec:Brownian_is_singular} we prove the main result of the paper, whilst in Section~\ref{sec:proof_Teo_aumann_integrall_null} we provide a proof to the statement \lq\lq zero--mean random set is a random element in $\Banach$ with zero--mean\rq\rq.

\section{Preliminaries}
\label{sec: preliminary results}

Here we refer mainly to \cite{li:ogu:kre02}. Denote by $\Closed[kc]{\Banach}$ the class of non--empty compact convex subsets of $\Banach$, endowed with the Hausdorff metric
\[
\HausDist(A,B) = \max\{\sup_{a\in A} \inf_{b\in B}\|a-b\| ,
\sup_{b\in B} \inf_{a\in A}\|a-b\|\},%
\]
and the operations
\[
A\Minkowski B = \{a+b:a\in A,\ b\in B\},\qquad  \lambda\cdot A=\lambda A =\{\lambda a: a\in A\}.
\]

A {\em fuzzy set} is a map $\nu: \Banach \to [0,1]$. Let $\Fuzzy[kc]{\Banach}$ denote the family of all fuzzy sets, which satisfy the following conditions.
\begin{enumerate}
\item
Each $\nu$ is an upper semicontinuous function, i.e.\ for each $\alpha\in (0,1]$, the cut set $\nu_\alpha=\{x\in\Banach : \nu (x)\ge \alpha \}$ is a closed subset of $\Banach$.

\item
The cut set $\nu_1=\{x\in\Banach : \nu (x)=1\}\neq\emptyset$.

\item
The support set $\nu_{0+}=\overline{\{x\in\Banach : \nu(x)>0\}}$ of $\nu$ is compact; hence every $\nu_\alpha$ is compact for $\alpha\in(0,1]$.

\item
For any $\alpha\in [0,1]$, $\nu_\alpha$ is a convex subset of $\Banach$.
\end{enumerate}
Let us endow $\Fuzzy[kc]{\Banach}$ with the metric
\[
\HausDist^{\infty}(\nu^1, \nu^2) = \sup \{\alpha\in [0,1] :  \HausDist (\nu^1_\alpha, \nu^2_\alpha)\}.
\]
and the operations
\[
(\nu^1\Mink\nu^2)_\alpha = \nu^1_\alpha \Minkowski \nu^2_\alpha,\qquad (\lambda\odot\nu^1)_\alpha = \lambda \cdot\nu^1_\alpha.
\]

Let $(\Omega,\salgebra,\prob)$ be a complete probability space. A {\em fuzzy set--valued random variable} (\frv{}) is a function $X:\Omega \to \Fuzzy[kc]{\Banach}$, such that $X_\alpha: \omega\mapsto X(\omega)_\alpha$ are random compact convex sets for every $\alpha\in (0,1]$ (i.e. $X_\alpha$ is a $\Closed[kc]{\Banach}$--valued function measurable with respect to the $\HausDist$--Borel $\sigma$--algebra).

An \frv{} $X$ is \emph{integrably bounded} and we shall write $X\in L^1[\Omega,\salgebra,\misura;\Fuzzy[kc]{\Banach}]=L^1[\Omega;\Fuzzy[kc]{\Banach}]$, if $\|X_{0+}\|_{H} :=\HausDist (X_{0+},\{0\})\in L^1[\Omega;\mathbb{R}]$.
\\
The \emph{expected value} of an \frv{} $X$, denoted by $\mathbb{E}[X]$, is a fuzzy set such that, for every $\alpha\in(0,1]$,
\[
(\mathbb{E}[X])_\alpha = \left(\int_\Omega X_\alpha \dd \misura \right)  =  \{\mathbb{E}(f): f\in L^1[\Omega; \Banach], f \in X_\alpha \ \mu-\textrm{a.e.}\}.
\]

\paragraph{Embedding Theorem.}

Let $\unitSphere$ be the unit sphere in $\Banach$. For any $\nu\in\Fuzzy[kc]{\Banach}$ define the \emph{support function} of $\nu$ as follows:
\[
\support[(x,\alpha)]{\nu} =\left\{
\begin{array}{ll}
\support[(x)]{\nu_\alpha} & {\rm if } \ \alpha >0, \\
\support[(x)]{\nu_{0^+}} & {\rm if } \ \alpha =0,
\end{array}
\right.
\]
for $(x,\alpha)\in\unitSphere\times [0,1]$ and where $\support[(x)]{K}=\sup\{ \langle x,a \rangle  : a\in K \}$, for  $x\in \unitSphere$.
\\
It is known that support function satisfies the following properties:
\begin{enumerate}
\item
for any $\nu^1, \nu^2 \in \Fuzzy[kc]{\Banach}$, $\support[(\cdot,\cdot)]{\nu^1\Mink \nu^2}=\support[(\cdot,\cdot)]{\nu^1}+\support[(\cdot,\cdot)]{\nu^2}$,

\item
for any $(x, \alpha)\in \Banach\times [0,1]$, $\support[(x,\alpha)]{X(\cdot)}\in L^1[\Omega;\mathbb{R}]$, $\mathbb{E}[\support[(x,\alpha)]{X}] = \support[(x,\alpha)]{\mathbb{E}[X]} $.
\end{enumerate}

Let $C(\unitSphere)$ denote the Banach space of all continuous functions $v$ on $\unitSphere$ with respect to the norm $\|v\|_C = \sup_{x\in\unitSphere} |v(x)|.$ Let $\overline{C}([0,1], C(\unitSphere) )$ be the set of all functions $f: [0,1]\to C(\unitSphere)$ such that $f$ is bounded, left continuous with respect to $\alpha\in (0,1]$, right continuous at 0, and $f$ has right limit for any $\alpha\in (0,1)$. Then we have that $\overline{C}([0,1], C(\unitSphere) )$ is a Banach space with the norm $\|f\|_{\overline{C}} = \sup_{\alpha\in [0,1]} \|f(\alpha) \|_C$, and the following embedding theorem holds.
\begin{pro}[(\cite{li:gua07} and the references therein.)]\label{pro:embedding}
There exists a function $j:\Fuzzy[kc]{\Banach} \to \overline{C}([0,1], C(\unitSphere) )$ such that:
\begin{enumerate}
\item
$j$ is an isometric mapping, i.e.
\[
\HausDist^\infty (\nu^1, \nu^2) = \| j(\nu^1) - j(\nu^2) \|_{\overline{C}}, \quad \nu^1, \nu^2\in\Fuzzy[kc]{\Banach},
\]

\item
$j(r\nu^1 + t\nu^2) = rj(\nu^1) + tj(\nu^2)$, $\nu^1, \nu^2\in\Fuzzy[kc]{\Banach}$ and $r,t \ge 0$.

\item
$j(\Fuzzy[kc]{\Banach})$ is a closed subset in $\overline{C}( [0,1] , C(\unitSphere) )$.
\end{enumerate}
\end{pro}
As a matter of fact, we can define an injection $j:\Fuzzy[kc]{\Banach} \to \overline{C} ([0,1],C(\unitSphere))$ by $j(\nu)=\support{\nu}$, i.e.\ $j(\nu)(x,\alpha) = \support[(x,\alpha)]{\nu}$ for every $(x,\alpha)\in \unitSphere\times [0,1]$, and this mapping $j$ satisfies above theorem. For simplification, let $\overline{\mathbf{C}}:=\overline{C} ([0,1],C(\unitSphere))$.
\\
From Proposition~\ref{pro:embedding} it follows that every \frv{} $X$ can be regarded as a random element of $\overline{\mathbf{C}}$ by considering $j(X)=\support{X} : \Omega \to \overline{\mathbf{C}}$, where $\support{X}(\omega)= \support{X(\omega)}$.

\paragraph{Fuzzy set--valued Brownian motion.}

For the results in this subsection we refer to \cite{li:gua07} or we shall specify if otherwise.
\begin{defi}[\cite{pur:ral85}]
A \frv{} $X:\Omega \to \Fuzzy[kc]{\Banach} $ is \emph{Gaussian} if $\support{X}$ is a Gaussian random element of $\overline{\mathbf{C}}$.
\end{defi}
A random element $\support{X}$ taking values in $\overline{\mathbf{C}}$ is Gaussian if and only if, for any $n\in\mathbb{N}$ and $f_1, f_2, \ldots, f_n\in \overline{\mathbf{C}}^*$, the real vector--valued random variable $( f_1(\support{X}), f_2(\support{X}), \ldots, f_n(\support{X}) )$ is Gaussian, where $\overline{\mathbf{C}}^*$ is the conjugate space of $\overline{\mathbf{C}}$ (i.e.\ the set of all continuous linear functionals on $\overline{\mathbf{C}}$).
\\
It follows from the properties of $\support{X}$ and elements in $\overline{\mathbf{C}}^*$ that $X+Y$ is Gaussian if $X$ and $Y$ are Gaussian \frv{}. Also $\lambda X$ is Gaussian whenever $X$ is Gaussian and $\lambda\in\mathbb{R}$.

\begin{pro}[{\cite[Theorem~6.1.7]{li:ogu:kre02}}] \label{pro:gaussian_process}
A \frv{} $X$ is Gaussian if and only if  $X$ is representable in the form
\[
X=\mathbb{E}[X]\Mink \indicator{\xi},
\]
where $\xi$ is a Gaussian random element of $\Banach$ with zero mean.
\end{pro}
\begin{defi}
Assume that $\{\salgebra_t: t\ge 0\}$ is a $\sigma$--filtration satisfying the usual condition (complete and right continuous). $\{X_t: t\ge 0\}$ is called an adaptive fuzzy set--valued stochastic process if for any $t\in\mathbb{R}_+$, $X_t$ is an $\salgebra_t$--measurable \frv{}. An adaptive fuzzy set--valued stochastic process $\{X_t : t\ge 0\}$ is called Gaussian if, for any $t\in\mathbb{R}_+$, $X_t$ is Gaussian.
\end{defi}
An adaptive fuzzy set--valued stochastic process $X=\{X_t : t\ge 0\}$ is Gaussian if and only if $\{ ( f_1(\support{X_t}), \ldots ,f_n (\support{X_t}) ) : t\ge 0 \}$ is a real vector--valued Gaussian process, for any $n\in\mathbb{N}$ and $f_1, f_2, \ldots, f_n \in \overline{\mathbf{C}}^*$. Further, the following theorem holds.
\begin{defi}
\label{def:B_t}
An adaptive fuzzy set--valued stochastic process $\{B_t: t\in\mathbb{R}_+\}$ is called a fuzzy set--valued Brownian motion if and only if $\{\support{B_t} :  t\in\mathbb{R}_+\}$ is a Brownian motion in $\overline{\mathbf{C}}$.
\end{defi}
\begin{pro}
\label{pro:caratterizzazione_Brownian}
Assume that a fuzzy set--valued stochastic process $\{B_t: t\ge 0\}$ satisfies $B_0=\indicator{0}$. Then $\{B_t : t\ge 0\}$ is a fuzzy set--valued Brownian motion if and only if it is a Gaussian process and
\begin{enumerate}
\item
$\mathbb{E}[f_i(\support{B_t})]=0$, for any $t\ge 0$, $f_i\in \overline{\mathbf{C}}^*$, $i=1,\ldots,n$,

\item
$\mathbb{E}[f_i(\support{B_t})f_i(\support{B_s})]=t\wedge s$, for any $s,t\ge 0$, $f_i\in \overline{\mathbf{C}}^*$, $i=1,\ldots,n$,

\item
$\mathbb{E}[f_i(\support{B_t})f_j(\support{B_s})]=0$, for any $s,t\ge 0$, $f_i,f_j\in \overline{\mathbf{C}}^*$, $i\neq j$, $i,j=1,\ldots,n$.

\end{enumerate}
\end{pro}

In \cite[Theorem~4.3 and Theorem~4.4]{li:gua07} the authors provide also some properties of a fuzzy set--valued Brownian motion that are very similar to those of the real case.
\begin{pro}
\label{pro:Brownian_properties}
Let $\{B_t : t\ge 0\}$ be a fuzzy set--valued Brownian motion. The following hold.
\begin{enumerate}
\item
$\{B_{t+t_0} \}_{t\ge 0}$ is a fuzzy set--valued Brownian motion for any $t_0\ge 0$.

\item\label{teo:2.Brownian_properties}
$\{\nu \Mink B_t \}_{t\ge 0}$ is a fuzzy set-valued Brownian motion for any fuzzy set $\nu\in\Fuzzy[k]{\Banach}$.

\item
$\{\frac{1}{\sqrt{\lambda}} B_{\lambda t} \}_{t\ge 0}$ is a fuzzy set-valued Brownian motion for any $\lambda> 0$.

\item
$\{t B_{\frac{1}{\sqrt{t}}} \}_{t\ge 0}$ is a fuzzy set-valued Brownian motion.

\item
If $\salgebra_t=\sigma\{B_s : s\le t\}$, then $\{B_t, \salgebra_t\}_{t\ge 0}$ is a fuzzy set--valued martingale.
\end{enumerate}
\end{pro}

\section[A \frv{} Brownian motion is Wiener in $\Banach$]{A \frv{} Brownian motion is a Wiener process in $\Banach$}%
\label{sec:Brownian_is_singular}

This section is devoted to prove Theorem~\ref{teo:Brownian_is_singular}: the main result of this paper.
\begin{teo}\label{teo:Brownian_is_singular}
A fuzzy set--valued process $\{B_t : t\ge 0\}$ is a Brownian motion, if and only if,
\[
B_t=\indicator{b_t}, \qquad \misura\textrm{--a.e.}
\]
where $\{b_t : t\ge 0\}$ is a Wiener process in $\Banach$.
\end{teo}
According to Definition~\ref{def:B_t} a fuzzy set--valued Brownian motion $B_t$ is a process taking values in $\Fuzzy{\Banach}$ (that is a functional space over $\Banach$). On the other hand, the previous result provides a way to handle a fuzzy set--valued Brownian motion simply using a random vector of $\Banach$. In other words, we observe a \lq\lq complexity reduction\rq\rq, i.e.\ from $\Fuzzy{\Banach}$ to $\Banach$.
\\
Moreover, in view of Theorem~\ref{teo:Brownian_is_singular}, Property~\ref{teo:2.Brownian_properties} in Proposition~\ref{pro:Brownian_properties} is true if and only if $\nu=\indicator{0}$,
whilst the remain properties in Proposition~\ref{pro:Brownian_properties} still hold due to the same properties of the driving Wiener process $b_t$ in $\mathbb{R}^d$.

Actually the \lq\lq complexity reduction\rq\rq\  stated in Theorem~\ref{teo:Brownian_is_singular} is strictly related to the characterization of Gaussian \frv{} (cf. Proposition~\ref{pro:gaussian_process}), to Property~1 of Proposition~\ref{pro:caratterizzazione_Brownian}, and to the following result obtained for random closed sets.
\begin{pro}\label{pro:aumann integral null}
Let $X$ be in $L^1[\Omega;\Closed{\Banach}]$ and let $a\in\Banach$. $\int_\Omega X d\misura=\{a\}$ if and only if there exists a $x\in L^1[\Omega;\Banach]$ such that $X=\{x\}$ $\misura$--a.e. and $\int_\Omega xd\misura=a$.
\end{pro}
\begin{cor}\label{cor: aumann integral null}
Let $X$ be in $L^1[\Omega;\Closed{\Banach}]$. $\int_\Omega X d\misura=\{0\}$ if and only if there exists a $x\in L^1[\Omega;\Banach]$ such that $X=\{x\}$ $\misura$-a.e. and $\int_\Omega xd\misura=0$.
\end{cor}
Although Proposition~\ref{pro:aumann integral null} and Corollary~\ref{cor: aumann integral null} are proved by Molchanov in \cite[Proposition~1.30, p.161]{mol05}, we shall propose in Appendix~\ref{sec:proof_Teo_aumann_integrall_null} alternative proofs via selections avoiding the use of the support function as Molchanov did.
\begin{lem}\label{lem:phi_in_dual}
For each $(x,\alpha)\in\Banach\times [0,1]$, the following map belongs to $\overline{\mathbf{C}}^*$
\[\begin{array}{rccl}
\varphi_{x,\alpha}:  & \overline{\mathbf{C}} & \to & \mathbb{R}\\
 & s & \mapsto & \varphi_{x,\alpha}(s)=s(x,\alpha).
\end{array}
\]
\end{lem}
\begin{myproof}
Map $\varphi_{x,\alpha}$ is linear since, for any $s_1$, $s_2$ in $\overline{\mathbf{C}}$ and $\lambda_1,\lambda_2\in\mathbb{R}$, the following chain of equalities hold.
\begin{align*}
\varphi_{x,\alpha}(\lambda_1 s_1+\lambda_2 s_2)= & [(\lambda_1 s_1+ \lambda_2 s_2)(\alpha)](x)= [\lambda_1 s_1(\alpha) + \lambda_2 s_2(\alpha)](x) \\
= & \lambda_1 s_1(\alpha,x) + \lambda_2 s_2(\alpha,x) = \lambda_1 \varphi_{x,\alpha}(s_1) + \lambda_2 \varphi_{x,\alpha}(s_2).
\end{align*}
For the continuity, let us consider any $s\in\overline{\mathbf{C}}$. For each $\varepsilon >0$ and $h\in\overline{\mathbf{C}}$ such that $\|h\|_{\overline{C}}<\varepsilon $, the following relations complete the proof.
\[
| \varphi_{x,\alpha}(s+h) - \varphi_{x,\alpha}(s) |= | \varphi_{x,\alpha}(h) | = | h(\alpha,x) | \le \|h\|_{\overline{C}} < \varepsilon.
\]
\end{myproof}

\begin{myproof}[ of Theorem~\ref{teo:Brownian_is_singular}]
The \lq\lq if\rq\rq  part is trivial.

In order to prove the \lq\lq only if\rq\rq  part let us consider the fuzzy set--valued Brownian motion $\{B_t : t\ge 0\}$.
\\
STEP 1. According to Proposition~\ref{pro:caratterizzazione_Brownian} and Proposition~\ref{pro:gaussian_process}, for any $t\ge 0$ and $f\in \overline{\mathbf{C}}^*$, it satisfies
\begin{align*}
0=\mathbb{E}[f(\support{B_t})]=\mathbb{E}[f(\support{\mathbb{E}[B_t]\Mink \indicator{\xi_t}})].
\end{align*}
where $\xi_t$ is an Gaussian random element of $\Banach$ with $\mathbb{E}\xi_t=0$. By the fact that, for any $\nu^1, \nu^2\in\Fuzzy[c]{\Banach}$, $\support{\nu^1\Mink \nu^2}= \support{\nu^1} + \support{\nu^2}$ (cf. Proposition~\ref{pro:embedding}), using the linearity of the expected value and of $f$, we get
\begin{align}
0 &= \mathbb{E}[f(\support{\mathbb{E}[B_t]})] + \mathbb{E}[f(\support{\indicator{\xi_t}})] = f(\support{\mathbb{E}[B_t]}) + f (\mathbb{E}[\support{\indicator{\xi_t}}]) \nonumber
\\
&= f(\support{\mathbb{E}[B_t]}) + f (\support{\indicator{\mathbb{E}[\xi_t]}}) = f(\support{\mathbb{E}[B_t]}), \label{eq:f(s)=0}
\end{align}
for any $t\ge 0$ and $f\in \overline{\mathbf{C}}^*$, where for the last two equalities we use $\support{\mathbb{E}{X}} = \mathbb{E} \support{X}$ and the fact that $\xi_t$ is zero mean.
\\
Clearly $\support{\mathbb{E}[B_t]}\equiv 0$. On the contrary, there will exists an $\alpha\in [0,1]$ such that $\support{\mathbb{E}[B_t]}(\alpha)\not\equiv 0 $; i.e.\ there exists an $\alpha\in [0,1]$ and $x\in\Banach$ such that $\support{\mathbb{E}[B_t]}(\alpha,x)\neq 0$. Let us consider the map defined by $\varphi_{x,\alpha}(s)=s(x,\alpha)$. It is an element of $\overline{\mathbf{C}}^*$ (cf. Lemma~\ref{lem:phi_in_dual}). Then $\varphi_{x,\alpha}(\support{\mathbb{E}[B_t]})\neq 0$ contradicts Equation~\eqref{eq:f(s)=0}.
\\%
As a consequence, $\mathbb{E}[B_t]=\indicator{0}$ for each $t\ge 0$; i.e.\ \begin{equation}\label{eq:Bt_alpha_level_is_null}
\mathbb{E}[(B_t)_\alpha]=\{0\},
\end{equation}
for each $t\ge 0$ and $\alpha\in (0,1]$.
\\%
STEP 2. Combining Corollary~\ref{cor: aumann integral null} with Equation~\eqref{eq:Bt_alpha_level_is_null} we obtain that, for each $t\ge 0$ and $\alpha\in (0,1]$, $(B_t)_\alpha$ is actually $\misura$--a.e. a random singleton with null mean value; i.e.\ $(B_t)_\alpha=\{b_t\}$ $\misura$--a.e. with $b_t$ being a random element of $\Banach$ such that $\mathbb{E}b_t=0$. By definition of $\alpha$--level sets for fuzzy set, $(B_t)_\alpha\supset (B_t)_\beta$ for any $0\le \alpha\le \beta\le 1$, and then $B_t=\indicator{b_t}$ $\misura$--a.e.\@.
\\%
Since $\{B_t\}_{t\ge 0}$ is a fuzzy set--valued Brownian motion, $\{b_t\}_{t\ge 0}$ is a Brownian motion in $\Banach$, and this fact concludes the proof.
\end{myproof}
Note that Proof of Theorem~\ref{teo:Brownian_is_singular} only uses the fact that $\{B_t\}$ is a Gaussian process for which any finite distribution, at any time $t$, has null expectation.

We want to point out that, although one can associate a fuzzy set--valued Brownian motion at any Brownian motion in $\overline{\mathbf{C}}$ (using the embedding in Proposition~\ref{pro:embedding}), in general, the contrary is not possible. This is due to the embedding properties. In fact, $j(\Fuzzy[kc]{\Banach})$ is a proper subset of $\overline{C}([0,1], C(\unitSphere) )$.
\\
As a consequence, a Gaussian element in $\overline{C}([0,1], C(\unitSphere) )$ can assume different values (even \lq\lq negative\rq\rq), whilst this could not happen in $\Fuzzy[kc]{\Banach}$ since, the embedding $j$ could not carry back all the possible \lq\lq fluctuations\rq\rq of gaussian element.

In this view, a definition of fuzzy set--valued Brownian motion, that take care completely the complexity of the (fuzzy) set--valued framework, has to take into account the above arguments and must pay attention to the possibly degeneracy.

\section{Proof of Proposition~\ref{pro:aumann integral null}}
\label{sec:proof_Teo_aumann_integrall_null}

In \cite[Proposition~1.30, p.161]{mol05} Molchanov proposed a proof of Proposition~\ref{pro:aumann integral null}. It involves the support function of a set.
Here we propose a different approach, via random sets selections, that is interesting by itself, and that leads to the same result.

For the sake of generality, here we shall consider ${\mathfrak{X}}$ to be a separable Banach space with $\borel[{\mathfrak{X}}]$ its borel $\sigma$--algebra and $(\Omega, \salgebra)$ to be a measurable space endowed with a positive finite measure $\misura$ (till now ${\mathfrak{X}}$ was $\mathbb{R}^d$ and $\misura$ a probability measure).
\\
In order to prove Proposition~\ref{pro:aumann integral null} we need the following two lemmas. Roughly speaking, the former says that any non--null vector in ${\mathfrak{X}}$ can be separated from zero using a suitable countable family of elements of $\dBanach$. The second lemma says that, for any couple of different (on some set of positive measure) integrable random elements in ${\mathfrak{X}}$, there exists an element of $\dBanach$ that separates (on a set of positive measure) these two random elements of ${\mathfrak{X}}$.
\begin{lem}\label{lem: exist separators in teo.aumann integral null}
There exists $\{\phi_n\}_{n\in\mathbb{N}}\subset\dBanach$ such that whenever $x\in{\mathfrak{X}}\setminus \{0\}$ there exists $n\in\mathbb{N}$ for which $\phi_n(x)\neq 0$.
\end{lem}
\begin{myproof}
Let $\{x_n\}_{n\in\mathbb{N}}$ be a dense subset of ${\mathfrak{X}}$. As a consequence of the Hahn-Banach Theorem (cf. \cite[Corollary~II.3.14, p.~65]{dun:sch58}) there exists $\{\phi_n\}_{n\in\mathbb{N}}\subset\dBanach$ such that
$\phi_n(x_n)=\|x_n\|_{{\mathfrak{X}}}$ and $\|\phi_n\|_{\dBanach}= 1$
for all $n\in\mathbb{N}$. Then
\begin{equation}\label{eq: lem1 aumann integral null}
-\|{y}\|_{{\mathfrak{X}}} \le \phi_n(y) \le \|{y}\|_{{\mathfrak{X}}}, \qquad \forall y\in{\mathfrak{X}}\setminus\{0\}, \forall n\in\mathbb{N}.
\end{equation}
Let $x\in{\mathfrak{X}}\setminus \{0\}$ and $n \in\mathbb{N}$ such that $\|{x-x_n}\|_{{\mathfrak{X}}}\le \frac{\|{x_n}\|_{{\mathfrak{X}}}}{2}$. By \eqref{eq: lem1 aumann integral null} we have
\[
\phi_n(x)=\phi_n(x_n) + \phi_n(x-x_n) \ge\|{x_n}\|_{{\mathfrak{X}}} - \|{x-x_n}\|_{{\mathfrak{X}}} \ge
\frac{\|{x_n}\|_{{\mathfrak{X}}}}{2}
> 0
\]
i.e.\ $\phi_n(x)>0$ that concludes the proof.
\end{myproof}
\begin{lem}\label{lem: varphi exists in teo.aumann integral null}
Let $x_1$, $x_2\in L^1[\Omega;{\mathfrak{X}}]$ and $A=\{\omega\in\Omega: x_1(\omega)\neq x_2(\omega)\}$ with $\misura(A) >0$. Then there exists $\varphi\in\dBanach$ such that
\[
A_\varphi = \{\omega\in \Omega: \varphi[x_1(\omega)] >
\varphi[x_2(\omega)]\}
\]
has positive measure (i.e.\ $\misura(A_\varphi)>0$).
\end{lem}
\begin{myproof}
Let $x=(x_1-x_2)$ then $A=\{\omega\in\Omega: x(\omega)\neq 0\}$ and let $\{\phi_n\}_{n\in\mathbb{N}}\subset\dBanach$ as in Lemma~\ref{lem: exist separators in teo.aumann integral null}. We claim that there exists $n\in\mathbb{N}$ such that $\misura (A_{\phi_n})+ \misura (A_{-\phi_n}) >0$.
By contradiction, if $A_n=A_{\phi_n}\cup A_{-\phi_n}$, we have
$$
\misura(A_n)\le \misura(A_{\phi_n}) + \misura(A_{-\phi_n}) = 0, \qquad
\forall n\in\mathbb{N}.
$$
Now we prove that $A\subseteq\bigcup_{n\in\mathbb{N}} A_n$: let $\omega\in A$ then $x(\omega)\neq 0$ and, by hypothesis, there exists $n\in\mathbb{N}$ such that $\phi_n(x(\omega))\neq 0$. Hence $\phi_n(x(\omega)) > 0$ or $\phi_n(x(\omega)) < 0$ i.e.\ $\omega\in A_n$ and thus $A\subseteq\bigcup_{n\in\mathbb{N}} A_n$.
\\%
This means that $\misura(A)\le\misura(\bigcup_{n\in\mathbb{N}} A_n)= 0$ that contradicts hypothesis ($\misura(A)>0$) and concludes the proof.
\end{myproof}
\begin{myproof}[ of Proposition~\ref{pro:aumann integral null}]
The ``if'' part is trivial. Vice versa, let us suppose that $\int_\Omega x d\misura=a$ holds for all $x\in S_X$, where integral is in the Bochner sense. Let us recall that a Bochner integrable map is also Pettis integrable and by definition (see \cite{tal84,mus91}) we have
\begin{equation}\label{eq: teo.aumann integral null}
\int_\Omega \phi(x)d\misura=\phi(a),\qquad  \forall\phi\in\dBanach,\
\forall x\in S_X.
\end{equation}
Now, by contradiction, let us suppose that $x_1,x_2$ are distinct elements of $S_X$ i.e.\ $A=\{\omega\in\Omega: x_1(\omega)\neq x_2(\omega)\}$ has positive measure. Then, by Lemma~\ref{lem: varphi exists in teo.aumann integral null}, there exists $\varphi\in\dBanach$ such that $A_\varphi = \{\omega\in \Omega: \varphi[x_1(\omega)] > \varphi[x_2(\omega)]\}$ has positive measure. Let us consider $ x_{\varphi}= \indicator{A_\varphi}x_1 + \indicator{A_\varphi^C}x_2$.
Clearly $x_\varphi$ is a selection of $X$ (i.e.\ $x_{\varphi}\in S_X$), and
\begin{align*}
\int_\Omega \varphi(x_{\varphi})d\misura = & \int_{A_\varphi}
\varphi(x_1)d\misura + \int_{A_\varphi^C} \varphi(x_2)d\misura
\\
> & \int_{A_\varphi} \varphi(x_2)d\misura + \int_{A_\varphi^C}
\varphi(x_2)d\misura=\varphi(a)
\end{align*}
which contradicts Pettis integrability \eqref{eq: teo.aumann integral null}.
\end{myproof}

\section{Conclusion}

We proved that a fuzzy set--valued Brownian motion is actually a degenerated process. In particular, it can actually be handle by a wiener process in the understanding space. This simplification is due mainly both to the well--known Gaussian degeneracy and to the \lq\lq null\rq\rq expectation.
\\
Moreover, we provided an alternative proof to Proposition~\ref{pro:aumann integral null}: an integrable set--valued map, which integral is a singleton, is almost everywhere an integrable singleton--valued map

We think that used hypothesis can be relaxed in different ways in order to get generalizations. For example, the space $\Banach$ can be replaced with a more general one. In this case, the difficulty lies in the fact that one have to redefine fuzzy set--valued Brownian motion in the new space as well as to use a different embedding theorem.


\begin{thebibliography}{10}


\bibitem{dun:sch58}
N.~Dunford and J.~T. Schwartz.
\newblock {\em Linear Operators. {P}art {I}}.
\newblock Wiley Classics Library. John Wiley \& Sons Inc., New York, 1988.

\bibitem{li:gua07}
S.~Li and L.~Guan.
\newblock {Fuzzy set--valued Gaussian processes and Brownian motions}.
\newblock {\em Information Sciences}, 177:3251--3259, 2007.

\bibitem{li:ogu:kre02}
S.~Li, Y.~Ogura, and V.~Kreinovich.
\newblock {\em Limit Theorems and Applications of Set--Valued and Fuzzy
  Set--Valued Random Variables}.
\newblock Kluwer Academic Publishers Group, Dordrecht, 2002.

\bibitem{mol05}
I.~Molchanov.
\newblock {\em Theory of random sets}.
\newblock Springer. (2005)

\bibitem{mus91}
K.~Musial.
\newblock {Topics in the theory of Pettis integration.}
\newblock {\em Rend. Ist. Mat. Univ. Trieste}, 23:177--262, 1991.

\bibitem{pur:ral85}
M.~L. Puri and D.~A. Ralescu.
\newblock The concept of normality for fuzzy random variables.
\newblock {\em Ann. Probab.}, 13:1373--1379, 1985.

\bibitem{tal84}
M.~Talagrand.
\newblock {Pettis integral and measure theory.}
\newblock {\em Mem. Am. Math. Soc.}, 307:224 p., 1984.

\end{thebibliography}
\end{document}